\documentclass{article}%
\usepackage{amssymb}
\usepackage{amsmath}
\usepackage{amsfonts}
\usepackage{graphicx}%
\usepackage{enumitem} 
\usepackage{url} 
\usepackage{textcomp} 
\usepackage[margin=1.25in]{geometry}
\providecommand{\U}[1]{\protect\rule{.1in}{.1in}}

\begin{document}

\begin{center}
\textbf{DECAY OF MASS OF THE SOLUTION TO THE CAUCHY PROBLEM OF THE 
$P$-LAPLACIAN WITH ABSORPTION ON INFINITE GRAPHS}

Alan A.Tedeev. Department of Mathematics, Sergo Ordzhonikidze Russian State University for Geological Prospecting, Moscow, Russia.
\end{center}

\begin{center}
	e-mail: alan.tedeev2013@gmail.com, ORCID: 0009-0009-5065-4102
\end{center}

\textbf{Abstract}

We consider the Cauchy problem for the nonstationary discrete p-Laplacian with absorption term on infinite graph which supports the Sobolev inequality. We established the precise conditions on parameters which guarantee the decay in time of the total mass
for nonnegative solutions. Our technique relies on suitable energy inequalities.

\textbf{Keywords}: infinite graphs, p-Laplacian,
absorption, decay mass, asymptotics for large times

Mathematics Subject Classification(2010) 35R02 \textperiodcentered\ 58J35
\textperiodcentered\ 39A12

\begin{center}
	1. INTRODUCTION
\end{center}

We consider nonnegative solutions to the Cauchy problem for discrete degenerate parabolic equations with absorption with potential of the form

\begin{equation}
\frac{\partial u}{\partial t}(x,t)=\Delta_{p}u(x,t)-q(x)u^{r}\text{, }x\in
V,\text{ }t>0,\text{ }(x,t)\in S_{T}=V\times(0,T)\text{,} \tag{1.1}%
\end{equation}

\begin{equation}
u(x,0)=u_{0}(x)\geq0\text{, }x\in V\text{.} \tag{1.2}%
\end{equation}
Here $V$ is the set of vertices of the graph $G(V,E,w)$ with edge set
$E\subset V\times V$ and weight $w$,%

\[
\Delta_{p}u(x,t):=\frac{1}{m(x)}%
{\displaystyle\sum\limits_{y\in Y}}
\left\vert u(y,t)-u(x,t)\right\vert ^{p-2}(u(y,t)-u(x,t))w(x,y)\text{, }%
\]

\[
q(x):=q(d(x,x_{0}))\text{, }x_{0}\in V\text{, }x_{0}\text{ be fixed,}%
\]
where $d(x,x_{0})$ is standard combinatorial distance in $G$ so that $d$ takes
only integral values. We assume that the graph $G$ is simple, undirected,
infinite, connected with locally finite degree%

\[
m(x):=%
{\displaystyle\sum\limits_{y\thicksim x}}
w(x,y)\text{,}%
\]
where we write $y\thicksim x$ if and only if $(x,y)\in E$. Here the weight
$w$: $V\times V\rightarrow\lbrack0,\infty)$ is symmetric, i.e.,
$w(x,y)=w(y,x)$, and strictly positive if and only if $y\thicksim x$, and
$w(x,x)=0,$ for $x\in V$. The density $q(x):V\rightarrow$ $\mathbb{R}%
\in\mathbb{N}$ is positive decreasing function which behavior at infinity will
be specified. In what follows we use the notation%

\[
D_{y}f(x)=f(y)-f(x)=-D_{x}f(y),\text{ }x,y\in V.
\]

We define :%

\[
\left\Vert f\right\Vert _{l^{q}(U)}^{q}=%
{\displaystyle\sum\limits_{x\in U}}
\left\vert f(x)\right\vert ^{q}m(x)\text{, }\left\Vert f\right\Vert
_{l^{\infty}(U)}=\sup_{x\in U}\left\vert f(x)\right\vert \text{, }\mu_{w}(U)=%
{\displaystyle\sum\limits_{x\in U}}
m(x)\text{.}%
\]

\textbf{Definition 1.1.} \textit{We say that }$G$\textit{\ satisfies a Sobolev
inequality for some given }$p\geq1$, $N>p$ \textit{if for any }$v>0$
\textit{and any finite subset} $U\subset V$\textit{\ with} $\mu_{w}%
(U)=v$\textit{\ here exists constant }$C$ \textit{independent of} $v$
\textit{such that} \textit{the following inequality holds true}%
\begin{equation}
\left(
{\displaystyle\sum\limits_{x\in U}}
\left\vert f(x)\right\vert ^{p^{\ast}}m(x)\right)  ^{1/p\ast}\leq C\left(
{\displaystyle\sum\limits_{x,y\in U}}
\left\vert f(y)-f(x)\right\vert ^{p}w(x,y)\right)  ^{1/p},\tag{1.3}%
\end{equation}
\textit{for all} \textit{real valued} $f:V\rightarrow R$ \textit{such that}
$f(x)=0$ if $x\notin U$. \textit{Here }$p^{\ast}=Np/(N-p)$.\textit{\ }

We assume additionally that for any $R\geq1$ there exists $C$ such that%

\begin{equation}
\mu(B(R)):=%
{\displaystyle\sum\limits_{x\in B(R)}}
m(x)\leq CR^{N}, \tag{1.4}%
\end{equation}
where $B(R)$ means a ball of radius $R$: $\left\{  x\in V:\text{ }d(x, x_{0})\leq
R\right\}  $, centered at fixed $x_{0}$.

Note that bound below of the volume growth $\mu(B(R))\geq CR^{N}$ follows (see for example [1-6]) from (1.3).

\textbf{Definition 1.2.} \textit{Let} $p>1$. \textit{We say that
}$u\in L^{\infty}(0,T;l^{1}(V))$\textit{\ \ and any }$T>0$ \textit{is a
solutions to (1.1), (1,2) if }$u(x,t)\in C^{1}([0,T])$\textit{\ for every
}$x\in V$\textit{\ and (1.1), (1.2) satisfy in classical point wise sense}.

Let us remark that a solution $u$ of the problem (1.1), (1.2) exists. For this
we refer the reader to [7,8]. In what follows we assume that%

\begin{equation}
2\leq p<r+1,\tag{1.5}%
\end{equation}
For some positive $0\leq\alpha_{1}\leq$ $\alpha_{2}<N$ we suppose that
$q(x)>0$ for $x\in V$ and for $s\geq1:$%

\begin{equation}
q(s)s^{\alpha_{2}}\text{ is nondecreasing, }q(s)s^{\alpha_{1}}\text{ is
nonincreasing} \tag{1.6}%
\end{equation}

In what follows we use a symbol $C$ means a generic constants which may vary from line
to line and depends on parameters of the problem $r,p,N,\alpha_{1},\alpha_{2}$
only. Further, when a function is defined only on the integers, we extend it to the positive real axis by linear interpolation. For brevity, we
use the same notation for both the original function and its extension. Denote
by $\Phi(R),$ $R>1$
\[
\Phi(R):=R^{p(r-1)/(r+1-p)}q(R)^{(p-2)/(r+1-p)}\text{,}%
\]

and the mass of solution by
\begin{equation*}
	M(t):=%
	{\displaystyle\sum\limits_{x\in V}}
	u(x,t)m(x).
\end{equation*}
If $\Phi(R)$ is increasing denote its inverse by $\widetilde{R}$. Our main
result reads as follows.

\textbf{Theorem}. \textit{Let $u(x,t)$ is the solution of (1.1) in $V\times(t>0)$, $u_{0}(x)\geq0$ and $u_{0}\in l_{1}
(V)$. Suppose that conditions (1.4) and (1.5) hold and $\Phi
(R)$ be an increasing function $R>1.$ Then for any $t>t_{0}$, for $t_{0}$ large enough, the following esimates hold}

\begin{equation}
M(t)\leq C\widetilde{R}(t)^{N}\widetilde{R}(t)^{-p/(r+1-p)}%
q(\widetilde{R}(t))^{-1/(r+1-p)}+%
{\displaystyle\sum\limits_{x\in V\diagdown B(\widetilde{R}(t))}}
u_{0}(x)m(x),\tag{1.6}%
\end{equation}

\begin{equation}
\left\Vert u(t)\right\Vert _{l^{\infty}(U)}\leq Ct^{-N/\lambda}M(t)^{p/\lambda
},\text{ }\lambda=N(p-2)+p\text{.} \tag{1.7}%
\end{equation}
\textit{where }$C$\textit{\ is independent of }$u_{0}(x)$.

\textbf{Corollary. }\textit{Provided}%

\[
\Phi(R)\rightarrow\infty\text{ as }R\rightarrow\infty\text{ and }%
C\widetilde{R}(t)^{N}\widetilde{R}(t)^{-p/(r+1-p)}q(\widetilde{R}%
(t))^{-1/(r+1-p)}\rightarrow0\text{ as }t\rightarrow\infty\text{,}%
\]
\textit{it follows (1.6) and (1.7) that}%

\[
M(t)\rightarrow0\text{ as }t\rightarrow\infty,
\]

\[
t^{N/\lambda}\left\Vert u(t)\right\Vert _{l^{\infty}(U)}\rightarrow0\text{ as
}t\rightarrow\infty.
\]

Consider particular cases of the Theorem.

1. Let $q=1.$ Then $\widetilde{R}(t)=t^{(r+1-p)/p(r-1)}$, $\widetilde
{R}(t)^{N}\widetilde{R}(t)^{-p/(r+1-p)}q(\widetilde{R}(t))^{-1/(r+1-p)}
\\= t^{\left(  N(r+1-p)-p\right)  /p(r-1)}\rightarrow0$ as
$t\rightarrow\infty$ if

\begin{equation}
r<r^{\ast}=p-1+\frac{p}{N}, \tag{1.8}%
\end{equation}

and thus $M(t)\rightarrow0$ as $t\rightarrow\infty.$

2. Let $q(R)=R^{-\alpha}$ and $H=p(r-1)-\alpha(p-2)$ be positive. Then $\widetilde{R}(t)=t^{(r+1-p)/H},$ 
and $M(t)\rightarrow0$ as $t\rightarrow\infty$ if%

\begin{equation}
r<r^{\ast}=p-1+\frac{p-\alpha}{N}\text{.} \tag{1.9}%
\end{equation}
In particular, if $\alpha=r$, then (1.9) reads%

\begin{equation}
r<r^{\ast}=\frac{N(p-1)+p}{N+1}\text{.}\tag{1.10}%
\end{equation}
Let us remark that in the continuous setting the critical exponents in (1.8),
(1.9) correspond to the Fujita critical exponent to the equation $u_{t}%
=\Delta_{p}u+u^{r}$ with blow-up term (see monograph [9] and references
therein). While, the critical exponents in (1.10) corresponds to the critical
exponents to the equation $u_{t}=\Delta_{p}u-\left\vert Du\right\vert ^{r}$, which is threshold between decay mass and no decay mass as $t\rightarrow
\infty$ (see [10,12]). In general, results of our theorem agree with
one obtained in [11,12].

\textbf{Example}. Let $G=\mathbb{Z}^{N},$ $q(R)=R^{-\alpha}\left[  \log R\right]
^{\beta},0\leq\alpha<p,$ $\beta>1$. Here, $\mathbb{Z}^{N}$ is a lattice graph. The set of vertices consists of all N-tuples $(x_{1},...x_{N})\sim(y_{1},...y_{N})$ if and only if 
\begin{equation*}
	\displaystyle\sum\limits_{i=1}^{N}\lvert x_{i}-y_{i}\rvert = 1.
\end{equation*}

Then 

$\Phi(R)=R^{H/(r+1-p)}\left[  \log
R\right]  ^{\beta(p-2)/(r+1-p)}$, and thus

\[
\widetilde{R}(s)=\frac{H}{r+1-p}s^{(r+1-p)/H}\left[  \log s\right]
^{-\beta(p-2)/H}(1+o(s))
\]

\[
C\widetilde{R}(t)^{N}\widetilde
{R}(t)^{-p/(r+1-p)}q(\widetilde{R}(t))^{-1/(r+1-p)}%
\]

\[
\leq Ct^{\left(  N(r+1-p)-p+\alpha\right)
/H}\left(  \log t\right)  ^{\gamma}\text{, }\gamma=\beta\left[  -(p-2)\frac
{N+\alpha-p}{H(r+1-p)}-\frac{1}{r+1-p}\right]  .
\]
\bigskip

In particular, if $\beta=0$, then

\[
C\widetilde{R}(t)^{N}\widetilde{R}(t)^{-p/(r+1-p)}q(\widetilde{R}(t))^{-1/(r+1-p)}
\]

\[
\leq Ct^{\left( N(r+1-p)-p+\alpha\right)/H}.
\]

To the best of our knowledge, our results are new even when $p=2.$

Due to the discrete nature of our problem, the well-known integral methods
commonly used in qualitative analysis in the continuous case (see for example
[11, 12] are not always applicable in our setting. This is due to the
specifics of the combinatorial distance. In particular, in the proofs of our statements, we avoid local iterations.

Let us remark also, that assumption (1.3) is essential to have $l^{1} \rightarrow l^{\infty}$ estimate as in (1.7), and we will not use it in the mass estimate (1.6).

Readers interested in the foundational aspects of combinatorial graph theory, particularly in relation to the elliptic and parabolic equations, may refer to the following works: [2–4]. While there is a substantial body of literature devoted to the elliptic problems on graphs—both linear and nonlinear—qualitative studies of nonlinear parabolic equations involving the $p$-Laplacian operator remain relatively scarce.

The analysis of the stationary $p$-Laplacian on infinite graphs has laid essential groundwork for understanding nonlinear dynamics in discrete settings. Foundational results by Holopainen and Soardi [13, 14] established Liouville-type theorems for $p$-harmonic functions, while Keller and Lenz [16] characterized spectral properties of unbounded Laplacians—insights critical for our study of the parabolic problem (1.1). Recent advances by Biagi, Meglioli, and Punzo [28] extended these results to elliptic equations with potentials under volume growth conditions (cf. (1.4)), and Ting and Feng [29] derived sharp $p$-Laplacian inequalities. Crucially, functional inequalities and embedding theorems, such as those by Pinchover-Tintarev [15] and Saloff-Coste [2, 3], underpin our Sobolev-type framework (Definition 1.1), enabling the energy methods used to prove mass decay. These works collectively highlight the interplay between graph geometry, potential terms, and nonlinear absorption, which we now adapt to the time-dependent case.

While blow-up phenomena for the $p$-Laplacian on finite graphs have been studied extensively (see, e.g., [19] for a survey), our work focuses on infinite graphs, where the interplay between diffusion, absorption, and graph geometry leads to qualitatively different behavior. For $p>2$, the mass conservation property was proven in [15], underscoring the critical role of nonlinearity in long-time dynamics. The linear case ($p=2$) introduces further complexity, as stochastic completeness and heat kernel estimates become central—questions addressed by Grigor'yan [4] and later refined by Andres-Deuschel-Slowik [22] for degenerate weights. Notably, Chung-Grigor'yan-Yau [6] linked spectral properties to isoperimetric inequalities, while Cardoso-Pinheiro [23] provided explicit spectral bounds for induced subgraphs, offering tools to quantify how local graph structure influences global behavior—a perspective relevant to our analysis of mass decay under condition (1.4).

For the Cauchy problem on infinite graphs, geometric conditions like Faber-Krahn inequalities (e.g., [7], [24]) have been used to study heat kernel decay, complementing the stationary $p$-Laplacian results in [25]. These works collectively highlight the delicate balance between graph structure (including spectral properties [23]), nonlinearity, and absorption that governs mass decay in (1.1).

For completeness, we mention that the discrete $p$-Laplacian on finite graphs has found applications in image processing through adaptive filtering techniques (Elmoataz et al. [20]) and data clustering [21]. Closer to our focus, the dynamics of solutions—including extinction rates and positivity preservation—were rigorously studied for finite graphs by Lee-Chung [17] and Xin et al. [18].

\begin{center}
	2. PROOF OF THEOREM 
\end{center}

We start with following technical lemma, which was proven in [24] :

\textbf{Lemma 2.1} \textit{Let }$\theta>0,$\textit{\ }$p>2$\textit{, }%
$h\geq0,u,v:$\textit{\ }$V\rightarrow R$\textit{. Then for all }$x,y\in V$%

\[
\mathit{(}\left\vert D_{y}u(x)\right\vert ^{p-2}\mathit{D}_{y}\mathit{u(x)-}%
\left\vert D_{y}v(x)\right\vert ^{p-2}\mathit{D}_{y}\mathit{v(x))}
\]

\begin{equation}
\mathit{\times D}_{y}\mathit{(u(x)-v(x)-h)}_{+}^{\theta}\mathit{\geq C}%
_{0}\left\vert D_{y}(u(x)-v(x)-h)_{+}^{(\theta-1+p)/p}\right\vert
^{p}\mathit{.} \tag{2.1}%
\end{equation}
We have%

\begin{equation}
M(t)=
{\displaystyle\sum\limits_{x\in B(R)}}
u(x,t)m(x)+
{\displaystyle\sum\limits_{x\in V\diagdown B(R)}}
u(x,t)m(x):=I_{1}(R)+I_{2}(R)\text{.} \tag{2.2}%
\end{equation}
Applying the H\"{o}lder inequality, we get%

\[
I_{1}(R)\leq\left(
{\displaystyle\sum\limits_{x\in B(R)}}
q(x)u(x,t)^{r}m(x)\right)  ^{1/r}\left(
{\displaystyle\sum\limits_{x\in B(R)}}
q(x)^{-1/(r-1)}m(x)\right)  ^{(r-1)/r}%
\]

\begin{equation}
\leq\left(  -\frac{dM(t)}{dt}\right)  ^{1/r}\varphi(R)^{(r-1)/r}%
\text{,}\tag{2.3}%
\end{equation}
where

\[
\varphi(R):=%
{\displaystyle\sum\limits_{x\in B(R)}}
q(x)^{-1/(r-1)}m(x)\text{.}%
\]
Thus, we have by (2.2) and (2.3):%

\begin{equation}
M(t)\leq\left(  -\frac{dM(t)}{dt}\right)  ^{1/r}\varphi(R)^{(r-1)/r}%
+I_{2}(R)\text{.}\tag{2.4}%
\end{equation}
Let $\zeta(d(x, x_{0}))$ be a cut-off function of a ball $B(R),$which we define as
follows $\zeta(d(x,x_{0}))=1$ when $d(x,x_{0})\geq3R$, $\zeta(d(x,x_{0}))=0$ for $d(x,x_{0})\leq2R$
and for $2R\leq d(x)\leq3R$:%

\[
\zeta=\frac{d(x,x_{0})-2R}{R}\text{.}%
\]
Multiplying both sides of the equation (1.1) by $\zeta^{s}(d(x,x_{0}))$,
$s>p$, and integrating by parts using Lemma 2.1, we get

\[%
{\displaystyle\sum\limits_{x\in V}}
u(x,t)\zeta^{s}m(x)+%
{\displaystyle\int\limits_{0}^{t}}
{\displaystyle\sum\limits_{x\in V}}
u^{r}q(x))\zeta^{s}m(x)d\tau
\]

\begin{equation}
=-\frac{s}{2}%
{\displaystyle\int\limits_{0}^{t}}
{\displaystyle\sum\limits_{x,y\in V}}
\zeta^{s-1}\left\vert D_{y}u\right\vert ^{p-2}D_{y}uD_{y}\zeta w(x,y)+%
{\displaystyle\sum\limits_{x\in V}}
u_{0}(x)\zeta^{s}m(x):=-\frac{s}{2}J_{1}+J_{2}\text{.} \tag{2.5}%
\end{equation}
Note that by symmetry of $w(x,y)$%

\[
J_{1}:=%
{\displaystyle\int\limits_{0}^{t}}
{\displaystyle\sum\limits_{x,y\in V}}
\zeta^{s-1}\left\vert D_{y}u\right\vert ^{p-2}D_{y}uD_{y}\zeta w(x,y)
\]%
\[
\leq%
{\displaystyle\int\limits_{0}^{t}}
{\displaystyle\sum\limits_{x,y\in V}}
\zeta^{s-1}(x)\left\vert D_{y}u\right\vert ^{p-1}\left\vert D_{y}%
\zeta\right\vert w(x,y)
\]%
\[
=%
{\displaystyle\int\limits_{0}^{t}}
{\displaystyle\sum\limits_{x,y\in V}}
\zeta^{s-1}(y)\left\vert D_{y}u\right\vert ^{p-1}\left\vert D_{y}%
\zeta\right\vert w(x,y).
\]
Thus,%

\[
J_{1}\leq%
{\displaystyle\int\limits_{0}^{t}}
{\displaystyle\sum\limits_{x,y\in V}}
\frac{1}{2}\left(  \zeta^{s-1}(y)+\zeta^{s-1}(x)\right)  \left\vert
D_{y}u\right\vert ^{p-1}\left\vert D_{y}\zeta\right\vert w(x,y):=\mathcal{E(}%
t,R)\text{,}%
\]
where the integrals are taken over $V_{R}=V\cap B(3R)\diagdown B(2R)$. By the H\"{o}lder inequality we get

\[
\mathcal{E(}t,R)\leq C\left(
{\displaystyle\int\limits_{0}^{t}}
\tau^{\beta}%
{\displaystyle\sum\limits_{x,y\in V_{R}}}
\zeta^{s}(x)\left\vert D_{y}u\right\vert ^{p}u^{-\sigma}w(x,y)d\tau\right)
^{(p-1)/p}%
\]

\begin{equation}
\times\left(  \frac{1}{R^{p}}%
{\displaystyle\int\limits_{0}^{t}}
\tau^{-\beta(p-1)}%
{\displaystyle\sum\limits_{x\in V_{R}}}
\zeta^{s-p}u^{\sigma(p-1)}m(x)d\tau\right)  ^{1/p}:=CJ_{3}^{(p-1)/p}%
R^{-1}J_{4}^{1/p}\text{,}\tag{2.6}%
\end{equation}
where we choose numbers $\beta$, $\sigma>0$: $\sigma(p-1)>1$, $\beta(p-1)<1.$

Let $\xi$ be a cut-off function such that $\xi=0$ if $d(x,x_{0})\leq R$ and
$d(x,x_{0})\geq4R$, $\xi=1$ if $\ 2R\leq d(x,x_{0})\leq3R$,%

\[
\xi=\frac{d(x,x_{0})-R}{R}\text{ if }R\leq d(x,x_{0})\leq2R\text{ and }%
\xi=\frac{4R-d(x,x_{0})}{R}\text{.}%
\]
Denote%

\[
I_{r+1-\sigma}:=%
{\displaystyle\int\limits_{0}^{t}}
\tau^{\beta}%
{\displaystyle\sum\limits_{x\in V}}
u^{r+1-\sigma}q(x)\xi^{s}m(x)d\tau\text{,}%
\]%
\[
J_{5}:=%
{\displaystyle\int\limits_{0}^{t}}
\tau^{\beta}%
{\displaystyle\sum\limits_{x,y\in V_{R}}}
\xi^{s}(x)\left\vert D_{y}u\right\vert ^{p}u^{-\sigma}w(x,y)d\tau.
\]
Next, multiplying both sides of (1.1) by $\tau^{\beta}u^{1-\sigma}\xi^{s}$ and
integrating by parts, we get%

\[
I_{r+1-\sigma}+J_{5}:=%
{\displaystyle\int\limits_{0}^{t}}
\tau^{\beta}%
{\displaystyle\sum\limits_{x\in V}}
u^{r+1-\sigma}q(x)\xi^{s}m(x)d\tau
\]

\[
+%
{\displaystyle\int\limits_{0}^{t}}
\tau^{\beta}%
{\displaystyle\sum\limits_{x,y\in V_{R}}}
\xi^{s}\left\vert D_{y}u\right\vert ^{p}u^{-\sigma}w(x,y)d\tau
\]

\[
\leq C%
{\displaystyle\int\limits_{0}^{t}}
\tau^{\beta-1}%
{\displaystyle\sum\limits_{x\in V}}
u^{2-\sigma}\xi^{s}m(x)d\tau
\]

\begin{equation}
+C%
{\displaystyle\int\limits_{0}^{t}}
\tau^{\beta}%
{\displaystyle\sum\limits_{x,y\in V}}
\xi^{s-1}\left\vert D_{y}u\right\vert ^{p-1}\left\vert D_{y}\xi\right\vert
u^{1-\sigma}w(x,y)d\tau:=CI_{2-\sigma}+CJ_{6}\text{.}\tag{2.7}%
\end{equation}
By the Young inequality we get%

\[
J_{6}\leq\frac{p-1}{p}\epsilon^{p/(p-1)}J_{5}+\frac{C}{R^{p}}\epsilon^{-p}%
{\displaystyle\int\limits_{0}^{t}}
\tau^{\beta}%
{\displaystyle\sum\limits_{x,y\in V}}
\xi^{s-p}u^{p-\sigma}\left\vert D_{y}\xi\right\vert ^{p}w(x,y)d\tau
\]%
\[
\leq\frac{p-1}{p}\epsilon^{p/(p-1)}J_{5}+\frac{C}{R^{p}}\epsilon^{-p}%
{\displaystyle\int\limits_{0}^{t}}
\tau^{\beta}%
{\displaystyle\sum\limits_{x,y\in V}}
\xi^{s-p}u^{p-\sigma}m(x)d\tau\text{.}%
\]%
\[
:=\frac{p-1}{p}\epsilon^{p/(p-1)}J_{5}+\frac{C}{R^{p}}\epsilon^{-p}%
I_{2-\sigma}.
\]
Absorbing the first term by small enough $\epsilon$, we get%

\[%
{\displaystyle\int\limits_{0}^{t}}
\tau^{\beta}%
{\displaystyle\sum\limits_{x\in V}}
u^{r+1-\sigma}q(x)\xi^{s}m(x)d\tau+J_{5}%
\]

\begin{equation}
\leq C%
{\displaystyle\int\limits_{0}^{t}}
\tau^{\beta-1}%
{\displaystyle\sum\limits_{x\in V}}
u^{2-\sigma}\xi^{s}m(x)d\tau+CI_{9}.\tag{2.8}%
\end{equation}

\[
I_{p-\sigma}=%
{\displaystyle\int\limits_{0}^{t}}
\tau^{\beta}%
{\displaystyle\sum\limits_{x,y\in V}}
\xi^{s-p}u^{p-\sigma}\left\vert D_{y}\xi\right\vert ^{p}w(x,y)d\tau
\]

\[
\leq\frac{C}{R^{p}}%
{\displaystyle\int\limits_{0}^{t}}
\tau^{\beta}%
{\displaystyle\sum\limits_{x\in V}}
\xi^{s-p}u^{p-\sigma}m(x)d\tau
\]

\begin{equation}
\leq\epsilon I_{r+1-\sigma}+\frac{C_{\epsilon}}{R^{p(r+1-\sigma)/(r+1-p)}}%
{\displaystyle\int\limits_{0}^{t}}
\tau^{\beta}%
{\displaystyle\sum\limits_{x\in B(4R)\diagdown B(R)}}
q(x)^{-(p-\sigma)/(r+1-p)}m(x)d\tau.\tag{2.9}%
\end{equation}
In what follows we choose $s=p+(p-\sigma)/(r+1-\sigma)>p$. In the same way we have%

\[
I_{2-\sigma}=%
{\displaystyle\int\limits_{0}^{t}}
\tau^{\beta-1}%
{\displaystyle\sum\limits_{x\in V}}
u^{2-\sigma}\xi^{s}m(x)d\tau
\]

\[
\leq\epsilon I_{r+1-\sigma}+C_{\epsilon}%
{\displaystyle\int\limits_{0}^{t}}
\tau^{\beta-(r+1-\sigma)/(r-1)}%
{\displaystyle\sum\limits_{x\in V}}
\xi^{s}q^{-(2-\sigma)/(r-1)}(x)m(x)d\tau\text{,}%
\]

\[
\frac{1}{R^{p}}%
{\displaystyle\int\limits_{0}^{t}}
\tau^{-\beta(p-1)}%
{\displaystyle\sum\limits_{x\in V}}
\xi^{s-p}u^{\sigma(p-1)}m(x)d\tau
\]

\begin{equation}
\leq\epsilon I_{r+1-\sigma}+C_{\epsilon}\frac{t^{1-\beta(p-1)(r+1)/(r+1-\sigma
p)}}{R^{p(r+1-\sigma)/(r+1-\sigma p)}}\left(
{\displaystyle\sum\limits_{x\in B(4R)\diagdown B(R)}}
q(x)^{-\sigma(p-1)/(r+1-\sigma p)}m(x)d\tau\right)  .\tag{2.10}%
\end{equation}
Summarizing (2.7)-(2.10), we have %

\[
I_{r+1-\sigma}\leq\frac{C}{R^{p(r+1-\sigma)/(r+1-p)}}%
{\displaystyle\int\limits_{0}^{t}}
\tau^{\beta}%
{\displaystyle\sum\limits_{x\in B(4R)\diagdown B(R)}}
q(x)^{-(p-\sigma)/(r+1-p)}m(x)d\tau
\]

\[
+Ct^{1+\beta-(r+1-\sigma)/(r-1)}%
{\displaystyle\sum\limits_{x\in V}}
\xi^{s}q^{-(2-\sigma)/(r-1)}(x)m(x)
\]

\[
\leq Ct^{1+\beta}R^{N}(R^{-p(r+1-\sigma)/(r+1-p)}q(R)^{-(p-\sigma
)/(r+1-p)}+t^{-(r+1-\sigma)/(r-1)}q(R)^{-(2-\sigma)/(r-1)})
\]

\[
=Ct^{1+\beta}R^{N}q(R)^{-(2-\sigma)/(r-1)}(R^{-p(r+1-\sigma)/(r+1-p)}%
q(R)^{-(p-\sigma)/(r+1-p)}%
\]

\[
+t^{-(r+1-\sigma)/(r-1)}q(R)^{-(2-\sigma)/(r-1)})^{r+1-\sigma}%
\]

\begin{equation}
=Ct^{1+\beta}R^{N}q(R)^{-(2-\sigma)/(r-1)}(\mathcal{A}(R,t))^{r+1-\sigma
},\tag{2.11}%
\end{equation}
where%

\[
\mathcal{A}(R,t):=R^{-p(r+1-\sigma)/(r+1-p)}q(R)^{-(p-\sigma)/(r+1-p)}%
+t^{-(r+1-\sigma)/(r-1)}q(R)^{-(2-\sigma)/(r-1)}.
\]
Here and hereafter we often use (1.4) and (1.6), which imply that for some $\theta < N$:
\begin{equation*}
	{\displaystyle\sum\limits_{x\in B(4R)\diagdown B(R)}}
	q(x)^{\theta}m(x) \leq C(\alpha_{1}, \alpha_{2})R^{N}q(R)^{\theta}.
\end{equation*}

Thus, applying the H\"{o}lder inequalty, we get

\[
\frac{1}{R^{p}}I_{p-\sigma}:=\frac{1}{R^{p}}%
{\displaystyle\int\limits_{0}^{t}}
\tau^{\beta}d\tau%
{\displaystyle\sum\limits_{x\in V}}
\xi^{s-p}u^{p-\sigma}m(x)
\]

\[
\leq\frac{1}{R^{p}}I_{r+1-\sigma}^{(p-\sigma)/(r+1-\sigma)}\left(
{\displaystyle\int\limits_{0}^{t}}
\tau^{\beta}d\tau%
{\displaystyle\sum\limits_{x\in B(4R)\diagdown B(R)}}
q(x)^{-(p-\sigma)/(r+1-p)}m(x)\right)  ^{(r+1-p)/(r+1-\sigma)}%
\]

\[
\leq C\frac{1}{R^{p}}\left(  t^{1+\beta}R^{N}q(R)^{-(2-\sigma)/(r-1)}\right)
^{(p-\sigma)/(r+1-\sigma)}\mathcal{A}(R,t)^{p-\sigma}%
\]

\[
\times t^{(1+\beta)(r+1-p)/(r+1-\sigma)}\left(  R^{N}q(R)^{-(p-\sigma
)/(r+1-\sigma)}\right)  ^{(r+1-p)/(r+1-\sigma)}%
\]

\begin{equation}
=Ct^{1+\beta}R^{N-p}q(R)^{-(p-\sigma)/(r-1)}\mathcal{A}(R,t)^{p-\sigma
},\tag{2.12}%
\end{equation}

\[
I_{2-\sigma}:=%
{\displaystyle\int\limits_{0}^{t}}
{\displaystyle\sum\limits_{x\in V}}
u^{2-\sigma}\tau^{\beta-1}\xi^{s}m(x)d\tau
\]

\[
\leq CI_{r+1-\sigma}^{(2-\sigma)/(r+1-\sigma)}\left(
{\displaystyle\int\limits_{0}^{t}}
\tau^{\beta-(2-\sigma)/(r-1)}q(R)^{-(2-\sigma)/(r-1)}d\tau%
{\displaystyle\sum\limits_{x\in x\in B(2R)\diagdown B(R)}}
m(x)\right)  ^{(r-1)/(r+1-\sigma)}%
\]

\[
\leq C\left(  t^{1+\beta}R^{N}q(R)^{-(2-\sigma)/(r-1)}(\mathcal{A}%
(R,t))^{r+1-\sigma}\right)  ^{(2-\sigma)/(r+1-\sigma)}%
\]

\[
\times\left(  t^{1+\beta-(2-\sigma)/(r-1)}\right)  ^{(r-1)/(r+1-\sigma
)}R^{N(r-1)/(r+1-\sigma)}%
\]

\begin{equation}
=Ct^{\beta}R^{N}q(R)^{-(2-\sigma)/(r-1)}\mathcal{A}(R,t))^{2-\sigma
},\tag{2.13}%
\end{equation}

\[
\frac{1}{R^{p}}I_{\sigma(p-1)}:=%
{\displaystyle\int\limits_{0}^{t}}
\tau^{-\beta(p-1)}%
{\displaystyle\sum\limits_{x\in V}}
\xi^{s-p}u^{\sigma(p-1)}m(x)d\tau\leq\frac{1}{R^{p}}I_{r+1-\sigma}%
^{\sigma(p-1)/(r+1-\sigma)}%
\]

\[
\times\left(
{\displaystyle\int\limits_{0}^{t}}
\tau^{-\beta(p-1)(r+1)/(r+1-\sigma p)}%
{\displaystyle\sum\limits_{x\in V}}
\xi^{s-p}q(R)^{-\sigma(p-1)/(r+1-\sigma p)}m(x)d\tau\right)  ^{(r+1-\sigma
p)/(r+1-\sigma)}%
\]

\[
=C\frac{1}{R^{p}}I_{r+1-\sigma}^{\sigma(p-1)/(r+1-\sigma)}\left(
t^{1-\beta(p-1)(r+1)/(r+1-\sigma p)}R^{N}q(R)^{-\sigma(p-1)/(r+1-\sigma
p)}\right)  ^{(r+1-\sigma p)/(r+1-\sigma)}%
\]

\[
\times\left(  t^{1+\beta}R^{N}q(R)^{-(2-\sigma)/(r-1)}(\mathcal{A}%
(R,t))^{r+1-\sigma}\right)  ^{\sigma(p-1)/(r+1-\sigma)}%
\]

\begin{equation}
=Ct^{1-\beta(p-1)}R^{N-p}\mathcal{A}^{\sigma(p-1)}q(R)^{-\sigma(p-1)/(r-1)}.
\tag{2.14}%
\end{equation}
Thus,%

\[
\mathcal{E(}t,R)\leq C(\frac{1}{R^{p}}%
{\displaystyle\int\limits_{0}^{t}}
\tau^{\beta}%
{\displaystyle\sum\limits_{x\in B(2R)\diagdown B(R)}}
\xi^{s}u^{p-\sigma}m(x)d\tau+
\]

\[%
{\displaystyle\int\limits_{0}^{t}}
\tau^{\beta-1}%
{\displaystyle\sum\limits_{x\in B(2R)\diagdown B(R)}}
\xi^{s}u^{2-\sigma}m(x)d\tau)^{(p-1)/p})%
\]

\[
\times\left(  \frac{1}{R^{p}}%
{\displaystyle\int\limits_{0}^{t}}
\tau^{-\beta(p-1)}%
{\displaystyle\sum\limits_{x\in V}}
\xi^{s-p}u^{\sigma(p-1)}m(x)d\tau\right)  ^{1/p}%
\]

\[
\leq C\left(  t^{\beta}R^{N}q(R)^{-\frac{2-\sigma}{r-1}}\mathcal{A}%
(R,t))^{2-\sigma}+t^{1+\beta}R^{N-p}q(R)^{-\frac{p-\sigma}{r-1}}%
\mathcal{A}(R,t)^{p-\sigma}\right)  ^{(p-1)/p}%
\]

\[
\times\left(  t^{1-\beta(p-1)}R^{N-p}\mathcal{A}^{\sigma(p-1)}q(R)^{-\sigma
(p-1)/(r-1)}\right)  ^{1/p}%
\]

\begin{equation}
=Ct^{1/p}R^{-1}R^{N}q(R)^{-2(p-1)/p(r-1)}\mathcal{A}^{2(p-1)/p}(1+tR^{-p}%
q(R)^{-(p-2)/(r-1)}\mathcal{A}^{p-2})^{(p-1)/p}\text{.}\tag{2.15}%
\end{equation}
Let $t>t_{0}>1$. Let $\widetilde{R}(t)$ be defined from %

\begin{equation}
R^{p(r-1)/(r+1-p)}q(R)^{(p-2)/(r+1-p)}=t.\tag{2.16}%
\end{equation}
Choose  $R=[\widetilde{R}(t)]$ - the  integer part of  $\widetilde{R}(t)$. Then%

\[
t\widetilde{R}(t)^{-p}q(\widetilde{R}(t))^{-(p-2)/(r-1)}\mathcal{A}%
^{p-2}(\widetilde{R}(t))\asymp1
\]
and we have%

\begin{equation}
\mathcal{E(}t,\widetilde{R}(t))\leq Ct^{1/p}\widetilde{R}^{-1}\widetilde
{R}^{N}q(\widetilde{R})^{-2(p-1)/p(r-1)}\mathcal{A}^{2(p-1)/p}(\widetilde
{R}).\tag{2.17}%
\end{equation}
Since, by definition%

\[
t^{1/p}\asymp\widetilde{R}^{(r-1)/(r+1-p)}(t)q(\widetilde{R}%
(t))^{(p-2)/p(r+1-p)}\text{,}%
\]
we have%

\begin{equation}
\mathcal{E(}t,\widetilde{R}(t))\leq C\widetilde{R}^{N}\widetilde
{R}^{-p/(r+1-p)}q(\widetilde{R})^{-1/(r+1-p)}\text{.}\tag{2.18}%
\end{equation}
Finally, since%

\[
\frac{1}{R^{p}}%
{\displaystyle\int\limits_{0}^{t}}
\tau^{-\beta(p-1)}%
{\displaystyle\sum\limits_{x\in V_{R}}}
\zeta^{s-p}u^{\sigma(p-1)}m(x)d\tau
\]

\[
\leq\frac{1}{R^{p}}%
{\displaystyle\int\limits_{0}^{t}}
\tau^{-\beta(p-1)}%
{\displaystyle\sum\limits_{x\in V_{R}}}
\xi^{s-p}u^{\sigma(p-1)}m(x)d\tau\text{,}%
\]
we conclude that%

\[%
{\displaystyle\sum\limits_{x\in V}}
\zeta^{s}u(x,t)m(x)\leq C\mathcal{E(}t,\widetilde{R}(t))+%
{\displaystyle\sum\limits_{x\in V}}
u_{0}(x)\zeta^{s}m(x)
\]

\begin{equation}
\leq C\widetilde{R}(t)^{N}\widetilde{R}(t)^{-p/(r+1-p)}q(\widetilde
{R}(t))^{-1/(r+1-p)}+%
{\displaystyle\sum\limits_{x\in V\diagdown B(\widetilde{R}(t))}}
u_{0}(x)m(x):=E(t).\tag{2.19}%
\end{equation}
Next, by (2.4) and the last estimate we get%

\[
M(\tau)\leq C\left(  -\frac{dM(\tau)}{d\tau}\right)  ^{1/r}\varphi
(R)^{(r-1)/r}+E(t)\text{, }0<\tau<t.
\]
Without loss of generality, we may assume that $M(\tau)\geq$ $E(t)$ for all
$\tau\in(0,t)$ otherwise nothing is to prove. Thus, we have%

\[
F(\tau):=M(\tau)-E(t)\leq C\left(  -\frac{dF(\tau)}{d\tau}\right)
^{1/r}\left(
{\displaystyle\sum\limits_{x\in B(\widetilde{R}(t))}}
q(x)^{-1/(r-1)}m(x)\right)  ^{(r-1)/r}.
\]
Integrating this inequality, we have%

\[
M(t)\leq E(t)+Ct^{-1/(r-1)}%
{\displaystyle\sum\limits_{x\in B(\widetilde{R}(t))}}
q(x)^{-1/(r-1)}m(x)
\]

\[
\leq E(t)+C\left(  t^{-1}\widetilde{R}(t)^{N(r-1)}q((\widetilde{R}%
(t))^{-1}\right)  ^{1/(r-1)}.
\]
Using also the definition of $\widetilde{R}(t)$, we conclude%

\[
M(t)\leq E(t)+C\widetilde{R}^{N}\widetilde{R}^{-p/(r+1-p)}q(\widetilde
{R})^{-1/(r+1-p)}%
\]

\[
\leq%
{\displaystyle\sum\limits_{x\in V\diagdown B(\widetilde{R}(t))}}
u_{0}(x)m(x)+C\widetilde{R}^{N}\widetilde{R}^{-p/(r+1-p)}q(\widetilde
{R})^{-1/(r+1-p)}.
\]
The first part of the Theorem is proved. To prove the second part we remark that we have [7], [24]

\[
\left\Vert u(t)\right\Vert _{\infty,V}\leq Ct^{-N/\lambda}M(t)^{p/\lambda
}\text{.}%
\]
Finally, using the mass estimate, we arrive at the desired result.

The proof of Theorem is complete. $\square$

\bigskip
\begin{center}
	REFERENCES
\end{center}

\begin{enumerate}
	\item \textit{Ostrovskii M.I.} Sobolev spaces on graphs // Quest. Math. 2005. Vol. 28. \textnumero 4. P. 501--523. DOI: 10.2989/16073600509486144.
	
	\item \textit{Saloff-Coste L.} Aspects of Sobolev-Type Inequalities. Cambridge: Cambridge University Press, 2002. (London Mathematical Society Lecture Note Series. Vol. 289).
	
	\item \textit{Saloff-Coste L., Pittet C.} Isophotometry and Volume Growth for Cayley Graphs. 2014. Available at: \url{https://pi.math.cornell.edu/lsc/papers/surv.pdf}.
	
	\item \textit{Grigor'yan A.} Analysis and Geometry on Graphs Part 1: Laplace operator on weighted graphs. Tsinghua University, 2012.
	
	\item \textit{Coulhon T., Grigor'yan A.} Random Walks on Graphs with Regular Volume Growth // Geom. Funct. Anal. 1998. Vol. 8. P. 656--701. DOI: 10.1007/s000390050070.
	
	\item \textit{Chung F., Grigor'yan A., Yau S.T.} Higher Eigenvalues and Isoperimetric Inequalities on Riemannian manifolds and graphs // Comm. Anal. Geom. 2000. Vol. 8. \textnumero 5. P. 969--1026. DOI: 10.1007/s000390050070.
	
	\item \textit{Mugnolo D.} Parabolic theory of the discrete p-Laplace operator // Nonlinear Anal. 2013. Vol. 87. P. 33--60. DOI: 10.1016/j.na.2013.04.002.
	
	\item \textit{Hua B., Mugnolo D.} Time regularity and long time behavior of parabolic p-Laplace equations on infinite graphs // J. Differ. Equ. 2015. Vol. 259. P. 6162--6190. DOI: 10.1016/j.jde.2015.07.018.
	
	\item \textit{Samarskii A.A., Galaktionov V.A., Kurdyumov S.P., Mikhailov A.P.} Blow-up in quasilinear parabolic equations. Berlin: de Gruyter, 1995. (de Gruyter Exp. Math. Vol. 19).
	
	\item \textit{Andreucci D., Tedeev A.F., Ughi M.} The Cauchy problem for degenerate parabolic equations with source and damping // Ukr Math Bull. 2004. Vol. 1. \textnumero 1. P. 1--23.
	
	\item \textit{Skrypnik I.I., Tedeev A.F.} Decay of the mass of the solution to the Cauchy problem of the degenerate parabolic equation with nonlinear potential // Complex Variables and Elliptic Equations. 2017. Vol. 63. \textnumero 1. P. 90--115. DOI: 10.1080/17476933.2017.1286331.
	
	\item \textit{Besaeva Z.V., Tedeev A.F.} The decay rate of the solution to the Cauchy problem for doubly nonlinear parabolic equation with absorption // Vladikavkaz Mathematical Journal. 2020. Vol. 22. \textnumero 1. P. 12--32. DOI 10.23671/VNC.2020.1.57535.
	
	\item \textit{Holopainen I., Soardi P.} A strong Liouville theorem for p-harmonic functions on graphs // Ann. Acad. Sci. Fenn. Math. 1996. Vol. 22.
	
	\item \textit{Holopainen I., Soardi P.M.} p-harmonic functions on graphs and manifolds // Manuscripta Math. 1997. Vol. 94. P. 95--110. DOI: 10.1007/BF02677841.
	
	\item \textit{Pinchover Y., Tintarev K.} On positive solutions of p-Laplacian-type equations // Analysis, Partial Differential Equations and Applications -- The Vladimir Maz'ya Anniversary Volume. Basel: Birkh\"auser, 2009. P. 225--236. (Operator Theory: Advances and Applications. Vol. 193).
	
	\item \textit{Keller M., Lenz D.} Unbounded Laplacians on graphs: basic spectral properties and the heat equation // Math. Model. Nat. Phenom. 2010. Vol. 5. P. 198--224.
	
	\item \textit{Xin Q., Mu C., Liu D.} Extinction and Positivity of the Solutions for a p-Laplacian Equation with Absorption on Graphs // J. Appl. Math. 2011. Art. ID 937079. DOI: 10.1155/2011/937079.
	
	\item \textit{Lee Y.S., Chung S.Y.} Extinction and positivity of solutions of the p-Laplacian evolution equation on networks // J. Math. Anal. Appl. 2012. Vol. 386. P. 581--592. DOI: 10.1016/j.jmaa.2011.08.023.
	
	\item \textit{Chung S.-Y., Choi M.-J.} Blow-up Solutions and Global Solutions to Discrete p-Laplacian Parabolic Equations // Abstr. Appl. Anal. 2014. Art. ID 351675. DOI: 10.1155/2014/351675.
	
	\item \textit{Elmoataz A., Toutain M., Tenbrinck D.} On the p-Laplacian and $\infty$-Laplacian on graphs with applications in image and data processing // SIAM J. Imaging Sci. 2015. Vol. 8. \textnumero 4. P. 2412--2451. DOI: 10.1137/15M1022793.
	
	\item \textit{Keller M., Lenz D., Schmidt M., Wirth M.} Diffusion determines the recurrent graph // Advances in Mathematics. 2015. Vol. 269. P. 364--398. DOI: 10.1016/j.aim.2014.10.003.
	
	\item \textit{Andres S., Deuschel J.D., Slowik M.} Heat kernel estimates for random walks with degenerate weights // Electron. J. Probab. 2016. Vol. 21. \textnumero 33. P. 1--21. DOI: 10.1214/16-EJP4382.
	
	\item \textit{Cardoso D.M., Pinheiro S.J.} Spectral Bounds for the k-Regular Induced Subgraph Problem // Appl. Comput. Matrix Anal. MAT-TRIAD 2015. Springer Proc. Math. Stat. 2017. Vol. 192. DOI: 10.1007/978-3-319-49984-0\_7.
	
	\item \textit{Andreucci D., Tedeev A.F.} Asymptotic estimates for the p-Laplacian on infinite graphs with decaying initial data // Potential Anal. 2020. Vol. 53. P. 677--699. DOI: 10.1007/s11118-019-09784-w.
	
	\item \textit{Keller M., Lenz D., Wojciechowski R.K.} Analysis and Geometry on Graphs and Manifolds. Cambridge: Cambridge University Press, 2020.
	
	\item \textit{Bianchi D., Setti G., Wojciechowski R.} The generalized porous medium equation on graphs: existence and uniqueness of solutions with l 1 data. 2022. arXiv: 2112.01733.
	
	\item \textit{Bianchi D., Setti A.G., Wojciechowski R.K.} The generalized porous medium equation on graphs: existence and uniqueness of solutions with l 1 data. 2022. arXiv: 2112.01733.
	
	\item \textit{Biagi S., Meglioli G., Punzo F.} A Liouville theorem for elliptic equations with a potential on infinite graphs // Calc. Var. 2024. Vol. 63. \textnumero 165. DOI: 10.1007/s00526-024-02768-8.
	
	\item \textit{Ting G.Y., Feng W.L.} p-Laplace elliptic inequalities on the graph // Commun. Pure Appl. Anal. 2025. P. 389--411. DOI: 10.3934/cpaa.2024094.
\end{enumerate}

\end{document}